\documentstyle[12pt]{article}
\def\vs#1{\noalign{\vskip#1pt}}
\def\dfrac{\displaystyle\frac}

\def\deg{\mathop{\rm deg}}
\def\UQ{U_q\bigl({\rm gl}(2/1)\bigr)}
\def\qqqq{\quad\qquad}
\def\ket#1{\bigl|#1\bigr>}

\def\be{\begin{equation}}
\def\ee{\end{equation}}

\begin{document}

\begin{center}
{\large\bf The $q$--boson--fermion realizations of the quantum
superalgebra $\UQ$}
\end{center}

\vskip0.6cm

\begin{center}
{\bf \v{C}. Burd\'{\i}k}

{\it Department of Mathematics and Doppler Institute, FNSPE,
Czech Technical University,

Trojanova 13, CZ-120 00 Prague 2, Czech Republic}

\vskip.3cm

{\bf O. Navr\'{a}til}

{\it Department of Mathematics, FTS,
Czech Technical University,

Na Florenci 25, CZ-110 00 Prague 1, Czech Republic}
\end{center}

\begin{abstract}
We show that our construction of realizations for Lie algebras and
quantum algebras can be generalized to quantum superalgebras,
too. We study an example of quantum superalgebra $\UQ$ and give
the boson--fermion realization with respect to one pair of
$q$--boson operators and 2 pairs of fermions.
\end{abstract}

\section{Introduction}

Boson--fermion realizations of a given set of operators via
Bose--Fermion creation and annihilation operators are among
the main tools of solving various quantum problems. The origin
is linked with the Schwinger \cite{Sch}, Dyson \cite{Dys} and
Holstein--Primakoff \cite{Hol} realizations which are different
boson realizations of the algebra ${\rm sl}(2)$.

Generalizations of the Dyson realization to the Lie algebra
${\rm sl}(n)$ were derived in \cite{O}. In our paper \cite{Bu1} we
formulated the method starting from the Verma modules for
obtaining boson realizations and in \cite{Bu2} we obtained
explicitly a braid class of realizations which generalized the
results from \cite{HaLas,EHa}.

Later the idea was extended to the Lie superalgebra, and the
Dyson type boson--fermion realizations were explicitly given in
\cite{Pa}, generalizing the results to ${\rm sl}(2/1)$
(\cite{AnLi},\cite{Ka}).

Today these boson--fermion realizations become a standard
technique in quantum many--body physics and we can also find
several other applications in all fields of quantum physics.

Quantum groups and quantum supergroups or $q$--deformed Lie
algebras and superalgebras imply some specific deformations
of the classical Lie algebras and superalgebras. From a
mathematical point of view, those are  noncommutative associative
Hopf algebras and superalgebras. The structure and representation
theory of quantum groups were extensively developed by Jimbo
\cite{Ji} and Drinfeld \cite{Dr}. The first "quantum" version of
Holstein--Primakoff was worked out for
$U_q\bigl({\rm sl}(2)\bigr)$ \cite{Chai} and then for
$U_q\bigl(({\rm sl}(3)\bigr)$ \cite{da}. The Schwinger type
realization was written in \cite{Mac} and \cite{Bie}. These
realizations found immediate applications [18--23].

In our papers \cite{Bu3,Bu4,Bu5} we studied the Dyson
realizations of the series algebras $U_q\bigl({\rm sl}(2)\bigr)$,
$U_q\bigl({\rm gl}(n)\bigr)$, $U_q(B_n)$, $U_q(C_n)$ and
$U_q(D_n)$. There is some special case \cite{Bu4} for which the
realization of the subalgebra $U_q\bigl({\rm gl}(n-1)\bigr)$ in
the recurrence is trivial. Such special realizations of the quantum
algebra $U_q\bigl({\rm sl}(n)\bigr)$ of Dyson type were studied in
\cite{Pa3}. 

The aim of the present paper is to show that there is a
possibility of generalizing our method \cite{Bu1} for deriving
the boson--fermion realization, too. This will be exemplified by
the quantum superalgebra $\UQ$. This superalgebra can be applied
to physical problems such as strongly correlated electron systems
\cite{bra1,bra2,KuSa}. We explicitly see the recurrence with
respect to $U_q\bigl({\rm gl}(1/1)\bigr)$ and consequently we
will show that again it is a generalization of the result from
\cite{Pa4}.

Some preliminary results concerning the general case
$U_q\bigl({\rm gl}(m/n)\bigr)$ have already been obtained and
prepared for publication.

\section{Preliminaries}

In this article, we will use the definition of a quantum
superalgebra $\UQ$ which can be found in \cite{Pa4}.

Let $q$ be an independent variable,
${\cal A}=C\bigl[q,q^{-1}\bigr]$ and ${\cal C}(q)$ be a division
field of ${\cal A}$. The superalgebra $\UQ$ is the associative
superalgebra over ${\cal C}(q)$ generated by even generators
$K_i$, $K_i^{-1}$, $i=1,\,2,\,3$, $E_{12}$, $E_{21}$ and odd
generators $E_{32}$, $E_{32}$ which satisfy the following
relations:
\be
\label{kr1}
\begin{array}{l}
K_i^{\pm1}K_j^{\pm1}=K_j^{\pm1}K_i^{\pm1}\,,\qquad
K_iK_i^{-1}=1\\
\vs{3}
K_iE_{jk}=q^{\delta_{ij}-\delta_{ik}}E_{jk}K_i\\
\vs{3}
\bigl[E_{12},E_{32}\bigr]=\bigl[E_{21},E_{23}\bigr]=0\\
\vs{3}
\bigl[E_{12},E_{21}\bigr]=
\dfrac{K_1K_2^{-1}-K_1^{-1}K_2}{q-q^{-1}}\\
\vs{3}
\bigl\{E_{23},E_{32}\bigr\}=
\dfrac{K_2K_3-K_2^{-1}K_3^{-1}}{q-q^{-1}}\\
\vs{3}
E_{23}^2=E_{32}^2=0\\
\vs{3}
E_{12}E_{13}-qE_{13}E_{12}=0\\
\vs{3}
E_{21}E_{31}-qE_{31}E_{21}=0\\
\end{array}
\ee
where
$$
\begin{array}{l}
E_{13}=E_{12}E_{23}-q^{-1}E_{23}E_{12}\\
\vs{3}
E_{31}=-E_{21}E_{32}+q^{-1}E_{32}E_{21}
\end{array}
$$
The Hopf structure of this superalgebra is defined by the
following operations:
\smallskip

\noindent
1. {\sl Coproduct\/} $\triangle$
$$
\begin{array}{ll}
\triangle(1)=1\otimes1\quad&
\triangle\bigl(K_i\bigr)=K_i\otimes K_i\\
\vs{3}
\triangle\bigl(E_{12}\bigr)=E_{12}\otimes K_1K_2^{-1}+
1\otimes E_{12}
\quad&
\triangle\bigl(E_{23}\bigr)=E_{23}\otimes K_2K_3+
1\otimes E_{23}\\
\vs{3}
\triangle\bigl(E_{21}\bigr)=E_{21}\otimes1+
K_1^{-1}K_2\otimes E_{21}
\quad&
\triangle\bigl(E_{32}\bigr)=E_{32}\otimes1+
K_2^{-1}K_3^{-1}\otimes E_{32}
\end{array}
$$

\noindent
2. {\sl Counit\/} $\varepsilon$
$$
\begin{array}{l}
\varepsilon(1)=\varepsilon\bigl(K_i\bigr)=1\\
\vs{3}
\varepsilon\bigl(E_{12}\bigr)=
\varepsilon\bigl(E_{23}\bigr)=
\varepsilon\bigl(E_{21}\bigr)=
\varepsilon\bigl(E_{32}\bigr)=0
\end{array}
$$

\noindent
3. {\sl Antipode\/} $S$
$$
\begin{array}{ll}
S(1)=1\quad&
S\bigl(K_i\bigr)=K_i^{-1}\\
\vs{3}
S\bigl(E_{12}\bigr)=-E_{12}K_1^{-1}K_2\quad&
S\bigl(E_{23}\bigr)=-E_{12}K_2^{-1}K_3^{-1}\\
\vs{3}
S\bigl(E_{21}\bigr)=-K_1K_2^{-1}E_{21}\quad&
S\bigl(E_{32}\bigr)=-K_2K_3E_{32}
\end{array}
$$
We do not use these operations for construction of the
realization.

The method of construction used is the same as in the case of the
Lie algebras \cite{Bu1} or quantum algebra \cite{Bu5} and is based
on using the induced representation. The difference from
quantum algebra is that together with $q$--deformed boson
operators \cite{Mac}, \cite{Bie} we also use fermion operators.

The algebra ${\cal H}$ of the $q$--deformed boson operators is
the associative algebra over the field ${\cal C}(q)$
generated by the elements of $a^+$, $a^{-}=a$, $q^x$ and
$q^{-x}$, satisfying the commutation relations
\be
\label{Hay1}
\begin{array}{lll}
q^xq^{-x}=q^{-x}q^x=1\,,\qquad&
q^xa^+q^{-x}=qa^+,\qquad&
q^xaq^{-x}=q^{-1}a,\\
aa^+-q^{-1}a^+a=q^x,\qquad&
aa^+-qa^+a=q^{-x},\quad&
\end{array}
\ee
The algebra ${\cal H}$ has faithful representation on vector
space with basic elements $\bigl\{\ket{n}$, where
$n=0,\,1,\,\dots\bigr\}$ of the form
\be
\label{Hay2}
q^x\ket{n}=q^n\ket{n}\,,\quad
a^+\ket{n}=\ket{n+1}\,,\quad
a\ket{n}=\bigl[n\bigr]\ket{n-1}\,,
\ee
where $\bigl[n\bigr]=\dfrac{q^n-q^{-n}}{q-q^{-1}}$.

Because of odd generators $E_{23}$ and $E_{32}$ we construct
realization by means of the algebra ${\cal H}$ for even
elements, and by fermion elements $b^+$ and $b$ for odd ones.
These fermion elements commute with the elements of ${\cal H}$
and together fulfil the relations
\be
\label{fer1}
bb=b^+b^+=0\,,\quad bb^++b^+b=1\,.
\ee

As in the case of the Lie algebras or quantum groups, our
realizations contain elements of quantum sub--superalgebra of
$\UQ$, namely, quantum superalgebra $U_q\bigl({\rm gl}(1/1)\bigr)$.
The element $x$ of this subalgebra commutes with the elements
from ${\cal H}$, and for the fermion elements $b^{\pm}$ the
relation
\be
\label{fer2}
xb^{\pm}=(-1)^{\deg x}b^{\pm}x\,,
\ee
holds.

Realization of the quantum superalgebra $\UQ$ is called the
homomorphism $\rho$ of the $\UQ$ to associative superalgebra
${\cal W}$ generated by ${\cal H}$, $b^{\pm}$ and
$U_q\bigl({\rm gl}(1/1)\bigr)$.

\section{Construction of the realization of $\UQ$}

First, for construction of the realization we find the induced
representation of $\UQ$. As subalgebra ${\cal A}_0$ of $\UQ$ we
choose a quantum superalgebra generated by $E_{23}$, $E_{21}$,
$E_{32}$, $K_i$ and $K_i^{-1}$, $i=1,\,2,\,3$. Let $\varphi$ be a
representation of ${\cal A}_0$ on vector space $V$. Let $\lambda$
be the left regular representation on $\UQ\otimes V$, i.e. for
$x,\,y\in\UQ$ and $v\in V$ the representation $\lambda$ is defined by
\be
\label{ind1}
\lambda(x)\bigl(y\otimes v\bigr)=xy\otimes v\,.
\ee
Let ${\cal I}$ be subspace of $\UQ\otimes V$ generated by the
relations
$$
xy\otimes v=x\otimes\varphi(y)v\,,
$$
for all $x\in\UQ$, $y\in{\cal A}_0$ and $v\in V$. It is easy to see
that the subspace ${\cal I}$ is $\lambda$--invariant. Therefore,
(\ref{ind1}) gives the representation on the factor--space
$W=\bigl[\UQ\otimes V\bigr]/{\cal I}$.

Let $E_{12}^NE_{13}^M=\ket{N,M}$. Due to the
Poincar\'e--Birkhoff-Witt theorem the space $W$ of the induced
representation is generated by the elements $\ket{N,M}\otimes v$
where $N=0,\,1,\,2,\,\dots$, $M=0,\,1$ and $v\in V$.

To obtain the explicit form of the induced representation, we
give some relations. They can be proved by mathematical induction
from relations (\ref{kr1}).

\noindent
{\bf Lemma 1.}
For any $n=0,\,1,\,2,\,\dots$ the following formulae hold:
$$
\begin{array}{l}
E_{13}E_{12}^n=q^{-n}E_{12}^nE_{13}\\
\vs{2}
E_{23}E_{12}^n=q^nE_{12}^nE_{23}-
q\bigl[n\bigr]E_{12}^{n-1}E_{13}\\
\vs{2}
E_{23}E_{13}^n=(-q)^nE_{13}^nE_{23}\\
\vs{2}
E_{32}E_{13}^n=(-1)^nE_{13}^nE_{32}+
\dfrac{1-(-1)^n}2\,q^{-n}E_{12}E_{13}^{n-1}K_2K_3\\
\vs{2}
E_{21}E_{12}^n=E_{12}^nE_{21}-
\dfrac{\bigl[n\bigr]}{q-q^{-1}}\,E_{12}^{n-1}
\bigl(q^{n-1}K_1K_2^{-1}-q^{-n+1}K_1^{-1}K_2\bigr)\\
\vs{2}
E_{21}E_{13}^n=E_{13}^nE_{21}+
\dfrac{1-(-1)^n}2\,E_{13}^{n-1}E_{23}K_1^{-1}K_2\\
\vs{2}
E_{31}E_{12}^n=E_{12}^nE_{31}+
q^{n-2}\bigl[n\bigr]E_{12}^{n-1}K_1K_2^{-1}E_{32}\\
\vs{2}
E_{31}E_{13}^n=(-1)^nE_{13}^nE_{31}+
\dfrac{1-(-1)^n}2\,q^{-1}E_{13}^{n-1}
\dfrac{K_1K_3-K_1^{-1}K_3^{-1}}{q-q^{-1}}\\
\vs{2}
E_{32}E_{23}^n=(-1)^nE_{23}^nE_{32}+
\dfrac{1-(-1)^n}2\,E_{23}^{n-1}
\dfrac{K_2K_3-K_2^{-1}K_3^{-1}}{q-q^{-1}}
\end{array}
$$

\medskip

We omit the details of the calculations and write the result for
the action of the induced representation on the basis elements
$\ket{N,M}\otimes v$.

\noindent
{\bf Theorem 1.}
The formulae
$$
\begin{array}{l}
E_{12}\ket{N,M}\otimes v=
\ket{N+1,M}\otimes v\\
\vs{3}
E_{13}\ket{N,M}\otimes v=
q^{-N_1}\ket{N,M+1}\otimes v\\
\vs{3}
E_{23}\ket{N,M}\otimes v=
-q\bigl[N\bigr]\,\ket{N-1,M+1}\otimes v+
(-1)^Mq^{N+M}\ket{N,M}\otimes\varphi\bigl(E_{23}\bigr)v\\
\vs{3}
K_1\ket{N,M}\otimes v=
q^{N+M}\ket{N,M}\otimes\varphi\bigl(K_1\bigr)v\\
\vs{3}
K_2\ket{N,M}\otimes v=
q^{-N}\ket{N,M}\otimes\varphi\bigl(K_2\bigr)v\\
\vs{3}
K_3\ket{N,M}\otimes v=
q^{-M}\ket{N,M}\otimes\varphi\bigl(K_3\bigr)v\\
\end{array}
$$
$$
\begin{array}{l}
E_{32}\ket{N,M}\otimes v=
\dfrac{1-(-1)^M}2\,q^{-M}\ket{N+1,M-1}\otimes
\varphi\bigl(K_2K_3\bigr)v+\\
\vs{3}
\qqqq+
(-1)^M\ket{N,M}\otimes\varphi(E_{32}\bigr)v\\
\vs{3}
E_{21}\ket{N,M}\otimes v=
-\dfrac{\bigl[N\bigr]q^{N+M-1}}{q-q^{-1}}\,\ket{N-1,M}\otimes
\varphi\bigl(K_1K_2^{-1}\bigr)v+\\
\vs{3}
\qqqq+
\dfrac{\bigl[N\bigr]q^{-N-M+1}}{q-q^{-1}}\,\ket{N-1,M}\otimes
\varphi\bigl(K_1^{-1}K_2\bigr)v+\\
\vs{3}
\qqqq+
\dfrac{1-(-1)^M}2\,\ket{N,M-1}\otimes
\varphi\bigl(E_{23}K_1^{-1}K_2\bigr)v+
\ket{N,M}\otimes\varphi\bigl(E_{21}\bigr)v\\
\vs{3}
E_{31}\ket{N,M}\otimes v=
\dfrac{1-(-1)^M}2\,q^{N-1}\bigl[N\bigr]\,\ket{N,M-1}\otimes
\varphi\bigl(K_1K_3\bigr)v+\\
\vs{3}
\qqqq+
(-1)^Mq^{N+M-2}\bigl[N\bigr]\,\ket{N-1,M}\otimes
\varphi\bigl(K_1K_2^{-1}E_{32}\bigr)v+\\
\vs{3}
\qqqq+
\dfrac{1-(-1)^M}2\,\dfrac{q^{-1}}{q-q^{-1}}\,\ket{N,M-1}\otimes
\Bigl(\varphi\bigl(K_1K_3-K_1^{-1}K_3^{-}\bigr)v+\\
\vs{3}
\qqqq+
(-1)^M\ket{N,M}\otimes\varphi\bigl(E_{31}\bigr)v
\end{array}
$$
give the induced representation of the quantum superalgebra
$\UQ$.

\medskip

We construct the realization of quantum superalgebra $\UQ$ from
the induced representation given in Theorem 1 as follows:

\noindent
We chose the representation $\varphi$, for which
$\varphi\bigl(E_{21}\bigr)v=0$, $\varphi\bigl(E_{31}\bigr)v=0$,
$\varphi\bigl(K_1\bigr)v=q^{\lambda}v$ and substitute
$$
\begin{array}{lll}
q^{\pm N}\to q^{\pm x}\quad&
\ket{N+1,M}\to a^+\quad&
\bigl[N\bigr]\,\ket{N-1,M}\to a\\
\vs{3}
\ket{N,M+1}\to b^+\quad&
\dfrac{1-(-1)^M}2\ket{N,M-1}\to b\quad&
q^{\pm M}\to\bigl(bb^++q^{\pm1}b^+b\bigr)\\
\vs{3}
\varphi\bigl(E_{21}\bigr)v\to0\quad&
\varphi\bigl(E_{31}\bigr)v\to0\quad&
\varphi\bigl(K_1^{\pm1}\bigr)v\to q^{\pm\lambda}\\
\vs{3}
\varphi\bigl(K_2^{\pm1}\bigr)v\to k_2^{\pm1}\quad&
\varphi\bigl(K_3^{\pm1}\bigr)v\to k_3^{\pm1}&\\
\vs{3}
(-1)^M\varphi\bigl(E_{23}\bigr)v\to e_{23}\quad&
(-1)^M\varphi\bigl(E_{32}\bigr)v\to e_{32}&
\end{array}
$$
(the last two relations reflect the fact that $e_{23}$ and
$e_{32}$ are fermions).

By this substitution we obtain the realization of the quantum
superalgebra $\UQ$.

\noindent
{\bf Theorem 2.}
The mapping $\rho:\UQ\to{\cal W}$ defined by the formulae
$$
\begin{array}{l}
\rho\bigl(E_{12}\bigr)=a^+\\
\vs{3}
\rho\bigl(E_{13}\bigr)=q^{-x}b^+\\
\vs{3}
\rho\bigl(E_{23}\bigr)=-qab^++q^x\bigl(bb^++qb^+b\bigr)e_{23}\\
\vs{3}
\rho\bigl(K_1\bigr)=q^{\lambda_1+x}\bigl(bb^++qb^+b\bigr)\\
\vs{3}
\rho\bigl(K_2\bigr)=q^{-x}k_2\\
\vs{3}
\rho\bigl(K_3\bigr)=\bigl(bb^++q^{-1}b^+b\bigr)k_3\\
\vs{3}
\rho\bigl(E_{32}\bigr)=q^{-1}a^+bk_2k_3+e_{32}\\
\vs{3}
\rho\bigl(E_{21}\bigr)=-\dfrac{a}{q-q^{-1}}
\Bigl(q^{\lambda_1+x-1}\bigl(bb^++qb^+b\bigr)k_2^{-1}-
q^{-\lambda_1-x+1}\bigl(bb^++q^{-1}b^+b\bigr)k_2\Bigr)-
q^{-\lambda_1}be_{23}k_2\\
\vs{3}
\rho\bigl(E_{31}\bigr)=a^+abq^{\lambda_1+x-1}k_3+
aq^{\lambda_1+x-2}\bigl(bb^++qb^+b\bigr)k_2^{-1}e_{32}+
q^{-1}b\dfrac{q^{\lambda_1}k_3-q^{-\lambda_1}k_3^{-1}}{q-q^{-1}}
\end{array}
$$
is the realization of the quantum superalgebra $\UQ$.
\smallskip

\noindent
This theorem can be proved by a direct calculation.

\section{Conclusion}

In this paper we gave the method of construction of the
$q$--boson--fermion realization of quantum superalgebras and
applied it to the quantum superalgebra $\UQ$. One of the
advantages of this method, in comparison with \cite{Pa4}, is that
we automatically obtain a realization and we do not need to
verify the generating relation. The reason is that the
representation of $q$--bosons and fermions on the vector space
$W$ with basis $\ket{N,M}$ is faithful.

The other advantage we see in the fact that our realization is
expressed by means of polynomials of $q$--deformed bosons and
fermions. On the other hand, we can easily obtain the Dyson
realization of quantum superalgebra. For this purpose, it is
sufficient to choose a realization of the generators of the
algebra ${\cal H}$ in the form
\be
\label{Hay3}
a^+=A^+\,,\quad
a=\dfrac{\bigl[N+1\bigr]}{N+1}\,A\,,\quad q^x=q^N\,,
\ee
where $\bigl[A,A^+\bigr]=1$ and $N=A^+A$. It is easy to verify
that the realization of $\UQ$ from Theorem 2 with realization
(\ref{Hay3}) of the algebra ${\cal H}$ and with a trivial
realization of subalgebra $U_q\bigl({\rm gl}(1/1)\bigr)$ leads,
after homomorphism of $\UQ$, to the realization given in
\cite{Pa4}. In this case, the realization is of course expressed
by means of a series in operators $A^+$ and $A$. Therefore,
we prefer our form of realizations.

Finally, our realizations contain, in contrast with those in
\cite{Pa4}, quantum sub-su\-per\-al\-geb\-ras. Various forms of
realizations of this sub-superalgebra give various realizations
of the quantum superalgebra. In the studied case, this
sub-superalgebra is $U_q\bigl({\rm gl}(1/1)\bigr)$, and, therefore,
is very simple. We can choose a realization of this superalgebra
as
$$
\rho\bigl(e_{23}\bigr)=\rho\bigl(e_{32}\bigr)=0\,,\quad
\rho\bigl(k_2\bigr)=\rho\bigl(k_3^{-1}\bigr)=q^{\lambda_2}\quad
{\rm and}\quad
\rho\bigl(k_2^{-1}\bigr)=\rho\bigl(k_3\bigr)=q^{-\lambda_2}\,.
$$
In this case, we obtain a realization with one $q$--deformed
boson pair, one fermion pair and two parameters.

However, by means of our method we construct other realization
of $U_q\bigl({\rm gl}(1/1)\bigr)$, namely, realization of the
form
$$
\begin{array}{l}
\rho\bigl(e_{23}\bigr)=b_2^+\\
\vs{3}
\rho\bigl(k_2\bigr)=q^{\lambda_2}\bigl(b_2b_2^++qb_2^+b_2\bigr)\\
\vs{3}
\rho\bigl(k_3\bigr)=q^{\lambda_3}\bigl(b_2b_2^++q^{-1}b_2^+b_2\bigr)\\
\vs{3}
\rho\bigl(e_{32}\bigr)=
\dfrac{q^{\lambda_2+\lambda_3}-q^{-\lambda_2-\lambda_3}}
{q-q^{-1}}\,b_2=
\bigl[\lambda_2+\lambda_3\bigr]b_2
\end{array}
$$
where $b_2$ and $b_2^+$ are the fermion elements.
If we use this realization of the quantum superalgebra in the
realization of $\UQ$ given in Theorem 2, we obtain realization
with one $q$--deformed boson pair, two fermion pairs and three
parameters, which corresponds to the case of the Lie and quantum
algebras.

As it is evident from \cite{Bu4,Bu5}, this method of
construction of realization is very successful for quantum groups.
Therefore, we believe that it will be very useful for
construction of realizations of quantum supergroups, too.

\bigskip

\noindent
{\small  Partial support from grant 201/01/0130
of the Czech Grant Agency is gratefully acknowledged.}

\end{document}